\documentclass[11pt,english]{amsart}
\usepackage[T1]{fontenc}
\usepackage[latin1]{inputenc}
\usepackage{babel}
\usepackage{amssymb}

\makeatletter

\providecommand{\LyX}{L\kern-.1667em\lower.25em\hbox{Y}\kern-.125emX\@}

 \theoremstyle{plain}    
 \newtheorem{thm}{Theorem} 
 \theoremstyle{plain}    
 \newtheorem{cor}{Corollary} 
 \theoremstyle{plain}    
 \newtheorem{lem}{Lemma} 
 \theoremstyle{plain}    
 \newtheorem{prop}{Proposition} 
 \theoremstyle{definition}
 \newtheorem{defn}{Definition}
 \theoremstyle{definition}
  \newtheorem{example}{Example}
 \theoremstyle{remark}
 \newtheorem{rem}{Remark}
 \theoremstyle{remark}    
 \newtheorem{claim}{Claim}
 \theoremstyle{remark}    
 \newtheorem*{acknowledgement*}{Acknowledgement} 

\usepackage{mathrsfs}
\let\mathcal\mathscr
\let\mathfrak\mathscr
\makeatother
\begin{document}

\title[classification of Poisson structures]{a classification of topologically stable Poisson structures on a
compact oriented surface}

\author{Olga Radko}

\begin{abstract}
Poisson structures vanishing linearly on a set of smooth closed disjoint
curves are generic in the set of all Poisson structures on a compact
connected oriented surface. We construct a complete set of invariants
classifying these structures up to an orientation-preserving Poisson
isomorphism. We show that there is a set of non-trivial infinitesimal
deformations which generate the second Poisson cohomology and such
that each of the deformations changes exactly one of the classifying
invariants. As an example, we consider Poisson structures on the sphere
which vanish linearly on a set of smooth closed disjoint curves. 
\end{abstract}
\maketitle

\section{Introduction}

Recently several results were obtained concerning the \emph{local
classification} of Poisson structures on a manifold. According to
the Splitting Theorem \cite{W-Local}, the problem of local classification
can be reduced to the classification of structures vanishing at a
point. In dimension \( 2 \), V. Arnold \cite{Arnold-book} obtained
a hierarchy of normal forms of germs of Poisson structures degenerate
at a point (see also a paper by  P.~Monnier \cite{Monnier} for a detailed exposition.)
Using the notion of the modular vector field of a Poisson structure,
\hbox{J.-P.}~Dufour and A.~Haraki \cite{Dufour-Haraki} and \hbox{Z.-J}.~Liu
and P.~Xu \cite{Liu-Xu-quadratic} obtained a complete local classification
of quadratic Poisson structures in dimension \( 3 \). Some results
related to local classification of Poisson structures in dimensions
\( 3 \) and \( 4 \) were also obtained by J.~Grabowski, G.~Marmo
and A.~M.~Perelomov in \cite{Grabowski_et_al}. 

However, not much is known in relation to the \emph{global classification}
of Poisson structures on a given manifold (i.e., classification up
to a Poisson isomorphism, see \cite{W-PG} for a general discussion).
In this paper we give an explicit example of a global classification
of a certain set of Poisson structures on a compact oriented surface
\( \Sigma  \). 

For a compact connected oriented surface \( \Sigma  \) and \( n\geq 0 \)
we consider the set \( \mathcal{G}_{n}(\Sigma ) \) of Poisson structures
on \( \Sigma  \) which vanish linearly on a set of \( n \) smooth
simple closed non-intersecting (in short, disjoint) curves. In particular,
\( \mathcal{G}_{0}(\Sigma ) \) is the space  of symplectic structures.
The set \( \mathcal{G}(\Sigma )\doteq \bigsqcup _{n\geq 0}\mathcal{G}_{n}(\Sigma ) \)
is dense in the vector space \( \Pi (\Sigma ) \) of all Poisson structures
on \( \Sigma  \). We call Poisson structures in \( \mathscr {G}(\Sigma ) \)
\emph{topologically stable} since the topology of their zero sets
is unchanged under small perturbations. Locally around a point on
a zero curve these structures have the first order of degeneracy in
the local Arnold classification.

The main result of this paper is a complete classification of the
topologically stable structures \( \mathcal{G}(\Sigma ) \) up to
an orientation-preserving Poisson isomorphism. We construct the set
of invariants of a structure \( \pi \in \mathcal{G}_{n}(\Sigma ) \)
which consists of

\begin{itemize}
\item the topological arrangement of the zero curves \( \gamma _{1},\dots ,\gamma _{n} \)
of \( \pi  \) taken with orientations defined by \( \pi  \) (up
to an orientation-preserving diffeomorphism);
\item the periods of the restriction of a modular vector field \( X^{\omega _{0}}(\pi ) \)
with respect to a volume form \( \omega _{0} \) to the zero curves
of \( \pi  \);
\item the volume invariant of \( \pi  \) generalizing the Liouville volume
for a symplectic structure.
\end{itemize}
We prove that two structures \( \pi ,\pi '\in \mathcal{G}_{n}(\Sigma ) \)
are globally equivalent via an orientation-preserving Poisson isomorphism
if and only if all the invariants for these structures coincide.

The question of global equivalence of Poisson structures by orientation-reversing
Poisson isomorphisms can be reduced to the orientation-preserving
case. Let \( \nu :\, \Sigma \to \Sigma  \) be an orientation-reversing
diffeomorphism of \( \Sigma  \). Two Poisson structures \( \pi ,\pi '\in \mathcal{G}_{n}(\Sigma ) \)
are globally equivalent via an orientation-reversing diffeomorphism
iff \( \nu _{*}\pi  \) and \( \pi ' \) are globally equivalent via
an orientation-preserving diffeomorphism.

We compute the Poisson cohomology of a given Poisson structure \( \pi \in \mathcal{G}_{n}(\Sigma ) \)
on a compact oriented surface \( \Sigma  \) of genus \( g \). The
zeroth cohomology (interpreted as the space of Casimir functions)
is generated by constant functions and is one-dimensional. The first
cohomology (interpreted as the space of Poisson vector fields modulo
Hamiltonian vector fields) has dimension \( 2g+n \) and is generated
by the image of the first de Rham cohomology of \( \Sigma  \) under
the injective homomorphism \( \tilde{\pi }:\, H^{*}(\Sigma )\to H^{*}_{\pi }(\Sigma ) \)
and by the following \( n \) vector fields:\[
X^{\omega _{0}}(\pi )\cdot \{\textrm{bump function around }\gamma _{i}\},\quad i=1,\ldots n. \]
(The map \( \tilde{\pi } \) on the level of cohomology is induced
by the canonical bundle map \( \tilde{\pi }:T^{*}\Sigma \to T\Sigma  \)
associated to the Poisson structure \( \pi  \) in the following way:
\( \beta (\tilde{\pi }(\alpha ))(x)\doteq \pi (\alpha ,\beta )(x) \),
where \( \alpha ,\beta \in \Gamma (T^{*}\Sigma ) \), \( x\in \Sigma  \).) 

The second cohomology is generated by a non-degenerate Poisson structure
\( \pi _{0} \) on \( \Sigma  \) and \( n \) Poisson structures
of the form\[
\pi _{i}=\pi \cdot \{\textrm{bump function around }\gamma _{i}\},\quad i=1,\dots ,n.\]
 Each of the generators of the second cohomology corresponds to a
one-parameter family of infinitesimal deformations of the Poisson
structure which affects exactly one of the numeric classifying invariants.
The deformation \( \pi \mapsto \pi +\varepsilon \cdot \pi _{0} \)
changes the regularized Liouville volume. For each \( i=1,\dots ,n \)
the deformation \( \pi \mapsto \pi +\varepsilon \cdot \pi _{i} \)
changes the modular period around the curve \( \gamma _{i} \). This
shows that the number of numeric classifying invariants \( (n+1) \)
for \( \mathcal{G}_{n}(\Sigma ) \)) equals to the dimension of the
second Poisson cohomology, and is, therefore, optimal. 

As an example we consider the topologically stable Poisson structures
\( \mathcal{G}(S^{2}) \) on the sphere. In this case, an explicit
description of the moduli space of topologically stable Poisson structures
up to an orientation-preserving Poisson isomorphism is obtained. 

Since topologically stable Poisson structures considered in the present
paper have degeneracies of the simplest kind, our work can be considered
as a first step in the direction of global classification of Poisson
structures having higher-order degeneracies in terms of the Arnold
hierarchy of germs of Poisson structures. We plan to pursue this direction
in our future work.

\begin{acknowledgement*}
I would like to thank my advisor, Prof. A.Weinstein, for proposing
this problem and for many fruitful discussions. 
\end{acknowledgement*}

\section{Topologically stable Poisson structures and their invariants}

First we recall the classification of symplectic structures on compact
oriented surfaces which follows from Moser's theorem.

\subsection{\label{sec:symplectic_structures}Classification of symplectic structures}

According to Darboux's theorem, all symplectic structures on a given
manifold \( M \) are locally equivalent: for each symplectic form
\( \omega  \) and a point \( p\in M^{2n} \), there exist a coordinate
system \( (U,\, x_{1},\cdots ,x_{n},\, y_{1},\cdots ,y_{n}) \) centered
at \( p \) such that \( \omega =\sum _{i=1}^{n}dx_{i}\wedge dy_{i} \)
on \( U \). Therefore, the dimension of the underlying  manifold is the only
local invariant of a symplectic structure. 

\begin{defn}
Two symplectic forms \( \omega  \) and \( \omega ' \) on
\( M \) are \emph{globally equivalent} if there is a symplectomorphism
\( \varphi :\, (M,\omega )\to (M,\omega ') \).
\end{defn}
In certain cases the following theorem of Moser allows one to classify
symplectic forms on a given manifold up to global equivalence:

\begin{thm}
(Moser, \cite{Moser}) Let \( \omega  \) and \( \omega ' \)
be symplectic forms on a compact connected manifold \( M \). Suppose
that \( [\omega ]=[\omega ']\in H^{2}(M) \) and that the \( 2 \)-form
\( \omega _{t}\doteq (1-t)\omega +t\omega ' \) is symplectic
for each \( t\in [0,1] \). Then there is a symplectomorphism \( \varphi :\, (M,\omega )\to (M,\omega ') \).
\end{thm}
The total Liouville volume \( V(\omega )\doteq \int _{M}\underbrace{\omega \wedge \cdots \wedge \omega }_{n} \)
associated to a symplectic structure \( \omega  \) on \( M \) is
a global invariant. That is, if symplectic forms \( \omega  \)
and \( \omega ' \) on a manifold \( M \) with \( \textrm{dim}(M)=2n \)
are globally equivalent, their Liouville volumes are equal, \( V(\omega )=V(\omega ') \).

In the case of a compact \( 2 \)-dimensional manifold Moser's theorem
implies

\begin{cor}
On a compact connected oriented surface \( \Sigma  \) two symplectic
structures \( \omega ,\, \omega '  \) are globally equivalent
iff the associated Liouville volumes are equal: \( \int _{\Sigma }\omega =\int _{\Sigma }\omega ' \).
\end{cor}

\subsection{\label{sec:description_of_structures}Description of the set \protect\( \mathfrak {G}(\Sigma )\protect \)
of topologically stable structures on \protect\( \Sigma \protect \)}

Let \( \Sigma  \) be a compact connected oriented \( 2 \)-dimensional
surface. Since there are no non-trivial \( 3 \)-vector fields, any bivector
field gives rise to a Poisson structure. Thus, Poisson structures
on \( \Sigma  \) form a vector space \( \Pi (\Sigma ) \). 

For \( n\geq 0 \) let \( \mathfrak {G}_{n}(\Sigma )\subset \Pi (\Sigma ) \)
be the set of Poisson structures \( \pi  \) on \( \Sigma  \) such
that

\begin{itemize}
\item the zero set \( Z(\pi )\doteq \{p\in \Sigma |\, \pi (p)=0\} \) of
\( \pi \in \mathfrak {G}_{n}(\Sigma ) \) consists of \( n \) smooth
disjoint curves \( \gamma _{1}(\pi ),\cdots ,\gamma _{n}(\pi ) \);
\item \( \pi  \) vanishes linearly on each of the curves \( \gamma _{1}(\pi ),\cdots ,\gamma _{n}(\pi ) \).
\end{itemize}
In particular, \( \mathfrak {G}_{0}(\Sigma ) \) is the set of symplectic
structures on \( \Sigma  \). Let \( \mathfrak {G}(\Sigma )\doteq \bigsqcup _{n\geq 0}\mathfrak {G}_{n}(\Sigma ) \).
The symplectic leaves of a Poisson structure \( \pi \in \mathfrak {G}(\Sigma ) \)
are the points in \( Z(\pi )=\bigsqcup ^{n}_{i=1}\gamma _{i} \) (the
\( 0 \)-dimensional leaves) and the connected components of \( \Sigma \setminus Z(\pi ) \)
(the \( 2 \)-dimensional leaves.) We call Poisson structures in \( \mathfrak {G} (\Sigma ) \) 
\emph {topologically stable} since the topology of their zero sets is unchanged under small perturbations.

Unless indicated otherwise, throughout the paper we denote by \( \omega _{0} \)
a symplectic form compatible with the orientation of \( \Sigma  \)
and by \( \pi _{0} \) the corresponding Poisson bivector. Since any
\( \pi  \) can be written as \( \pi =f\cdot \pi _{0} \) for a function
\( f\in C^{\infty }(\Sigma ) \), we have \( \Pi (\Sigma )=C^{\infty }(\Sigma )\cdot \pi _{0} \).
The subspace \( \mathfrak {G}_{n}(\Sigma ) \) corresponds in this
way to the product \( \mathfrak {F}_{n}(\Sigma )\cdot \pi _{0} \),
where \( \mathfrak {F}_{n}(\Sigma ) \) is the space of smooth functions
for which \( 0 \) is a regular value and whose zero set consists
of \( n \) smooth disjoint curves. 

According to the Elementary Transversality Theorem (see, e.g., Corollary 4.12 in \cite{GG}), the set  \( \mathfrak {F} (\Sigma )=\bigsqcup _{n\geq 0}\mathfrak {F}_{n}(\Sigma ) \) of functions for which \( 0\in \mathbb {R} \) is a regular value  forms an open dense subset of  \( C^{\infty }(\Sigma ) \) in the Whitney
\( C^{\infty } \) topology. Therefore, we have the following

\begin{prop}
The set of Poisson structures \( \mathfrak {G}(\Sigma ) \) is generic
inside of \( \Pi (\Sigma ) \), i.e. \( \mathfrak {G}(\Sigma ) \)
is an open dense subset of the space \( \Pi (\Sigma ) \) of all Poisson
structures on \( \Sigma  \) endowed with the Whitney \( C^{\infty } \)
topology.
\end{prop}
Of course, for some sets of disjoint curves on a surface there are
no functions (and, therefore, no Poisson structures) vanishing linearly
on that set and not zero elsewhere. For example, such is the case
of one non-separating curve on a \( 2 \)-torus.

We will use the following definition

\begin{defn}
Two Poisson structures \( \pi  \) and \( \pi ' \) on an oriented
manifold \( P \) are \emph{globally equivalent} if there is an orientation-preserving
Poisson isomorphism \( \varphi :\, (P,\pi )\to (P,\pi ') \).
\end{defn}
The main goal of this paper is to classify the set \( \mathcal{G}(\Sigma ) \)
of topologically stable Poisson structures up to global equivalence. 

First, we will need the following

\begin{lem}
\label{claim:orientation_curves} A topologically stable Poisson structure
\( \pi \in \mathcal{G}_{n}(\Sigma ) \) defines an orientation on
each of its zero curves \( \gamma _{i}\in Z(\pi ),\, i=1,\dots ,n \).
Moreover, this induced orientation on the zero curves of \( \pi  \)
does not depend on the choice of orientation of \( \Sigma  \).
\end{lem}
\begin{proof}
Let \( \omega _{0} \) be a symplectic form on \( \Sigma  \), and
\( \pi _{0}=(\tilde{\omega }_{0}^{-1}\otimes \tilde{\omega }_{0}^{-1})(\omega _{0}) \)
be the corresponding Poisson bivector. Since \( \pi =f\cdot \pi _{0} \)
and \( f \) vanishes linearly on each of \( \gamma _{i}\in Z(\pi ) \)
and nowhere else, \( f \) has constant sign on each of the \( 2 \)-dimensional
symplectic leaves of \( \pi  \). In particular, \( f \) has the
opposite signs on two leaves having a common bounding curve \( \gamma _{i} \).
This defines an orientation on \( \gamma _{i} \) in the following
way. For a non-vanishing vector field \( X \) tangent to the curve
\( \gamma _{i} \), we say that \( X \) is \emph{positive} if \( \omega _{0}(X,Y)\geq 0 \)
for all vector fields \( Y \) such that \( L_{Y}f\geq 0 \). We say
that \( X \) is \emph{negative} if \( -X \) is positive.

Suppose that \( X \) is a vector field tangent to \( \gamma _{i} \)
and positive on \( \gamma _{i} \) with respect to the chosen orientation
of \( \Sigma  \). If \( \omega _{0}' \) is a symplectic form inducing
the opposite orientation on \( \Sigma  \), then \( \omega '_{0}=-\alpha \cdot \omega _{0} \)
with \( \alpha \in C^{\infty }(\Sigma ), \) \( \alpha >0 \), and
\( \pi =-\alpha \cdot f\cdot \pi _{0} \). Since for \( Y' \) such
that \( L_{Y'}(-\alpha \cdot f)\geq 0 \) we have \( L_{Y'}f\leq 0 \),
it follows that \( \omega _{0}(X,Y')\leq 0 \) and, therefore,\[
\omega _{0}'(X,Y')=-\alpha \cdot \omega _{0}(X,Y')\geq 0.\]
Hence, if \( X \) is positive on \( \gamma _{i} \) with respect
to a chosen orientation of \( \Sigma  \), it is also positive on
\( \gamma _{i} \) with respect to the reverse orientation of \( \Sigma  \). 
\end{proof}
We will refer to this orientation of \( \gamma _{i}\in Z(\pi ) \)
as the \emph{orientation defined by \( \pi  \).}

\subsection{\label{sec:equivalence_curve_sets}Diffeomorphism equivalence of
sets of disjoint oriented curves}

We will use the following definition:

\begin{defn}
Two sets of smooth disjoint oriented curves \( (\gamma _{1},\dots ,\gamma _{n}) \)
and \( (\gamma _{1}',\dots ,\gamma _{n}') \) on an oriented surface
\( \Sigma  \) are called \emph{diffeomorphism equivalent} (denoted
by \( (\gamma _{1},\dots ,\gamma _{n})\sim (\gamma _{1},\dots ,\gamma _{n}') \))
if there is an orientation-preserving diffeomorphism \( \varphi :\, \Sigma \to \Sigma  \)
mapping the first set onto the second one and preserving the orientations
of curves. That is to say, for each \( i\in 1,\dots ,n, \) there exists
\( j\in 1,\dots ,n, \) such that \( \varphi (\gamma _{i})=\gamma '_{j} \)
(as oriented curves.) 
\end{defn}
Let \( \mathfrak {C}_{n}(\Sigma ) \) be the space of \( n \) disjoint
oriented curves on \( \Sigma  \) and \( \mathfrak {M}_{n}(\Sigma ) \)
be the moduli space of \( n \) disjoint oriented curves on \( \Sigma  \)
modulo the diffeomorphism equivalence relation, \( \mathfrak {M}_{n}(\Sigma )=\mathfrak {C}_{n}(\Sigma )/\sim  \).
For a set of disjoint oriented curves \( (\gamma _{1},\cdots ,\gamma _{n}) \),
let \( [(\gamma _{1},\cdots ,\gamma _{n})]\in \mathfrak {M}_{n}(\Sigma ) \)
denote its class in the moduli space \( \mathfrak {M}_{n}(\Sigma ) \).
If ~\( \bigsqcup _{i=1}^{n}\gamma _{i}=Z(\pi ) \) for a Poisson
structure \( \pi  \), we will also write \( [Z(\pi )]\in \mathfrak {M}_{n} \)
to denote the class of the set of curves \( (\gamma _{1},\dots ,\gamma _{n}) \)
taken with the orientations defined by \( \pi  \). 

The topology of the inclusion \( Z(\pi )\subset \Sigma  \) and the
orientations of the zero curves of a topologically stable Poisson
structure \( \pi \in \mathcal{G}(\Sigma ) \) are clearly invariant under
orientation-preserving Poisson isomorphisms. In other words, if \( \pi ,\pi '\in \mathcal{G}_{n}(\Sigma ) \)
are globally equivalent, \( [Z(\pi )]=[Z(\pi ')]\in \mathcal{M}_{n}(\Sigma ) \).

\subsection{\label{sec:modular}Modular period invariants}

Recall the definition and some properties of a  modular vector field
of a Poisson manifold, introduced by A.~Weinstein in \cite{Weinstein_modular}.

\begin{defn}
(Weinstein, \cite{Weinstein_modular}) For a volume form \( \Omega  \)
on an orientable Poisson manifold \( (P,\pi ) \) the \emph{modular
vector field} \( X^{\Omega } \) of \( (P,\pi ) \) with respect to
\( \Omega  \) is defined by\[
X^{\Omega }\cdot h\doteq \frac{L_{X_{h}}\Omega }{\Omega },\qquad h\in C^{\infty }(P).\]
Here \( X_{h}=\tilde{\pi }(dh) \) is the hamiltonian vector field
of \( h \), where \( \tilde{\pi }:\, T^{*}M\to TM \) is the canonical
bundle map associated to the Poisson bivector \( \pi  \). 
\end{defn}
Modular vector fields have the following properties (see \cite{Weinstein_modular}
for details):

\begin{enumerate}
\item \( X^{\Omega }=0 \) iff \( \Omega  \) is invariant under the flows
of all hamiltonian vector fields of \( \pi  \);
\item For any other \( \Omega ' \) the difference \( X^{\Omega '}-X^{\Omega } \)
is the hamiltonian vector field \( X_{\log \left| \Omega ' / {\Omega } \right|  } \);
\item The flow of \( X^{\Omega } \) preserves \( \pi  \) and \( \Omega  \), i.e., 
\( L_{X^{\Omega }}\pi =0 \) and  \( L_{X^{\Omega }}\Omega =0 \);
\item \( X^{\Omega } \) is tangent to the symplectic leaves of maximal
dimension.
\end{enumerate}
Let now \( \pi \in \mathfrak {G}_{n}(\Sigma ) \) be a topologically
stable Poisson structure on a surface \( \Sigma  \) as above. A symplectic
form \( \omega _{0} \) compatible with the orientation of \( \Sigma  \)
is also a volume form on \( \Sigma  \). Since the modular vector
field \( X^{\omega _{0}} \) preserves \( \pi  \), it follows that
the restriction of \( X^{\omega _{0}} \) to a curve \( \gamma _{i}\in Z(\pi ) \)
is tangent to \( \gamma _{i} \) for each \( i\in 1,\dots ,n \).
Since for another volume form \( \omega _{0}' \) the difference \( X^{\omega _{0}}-X^{\omega _{0}'} \)
is a hamiltonian vector field and, therefore, vanishes on the zero
set of \( \pi  \), it follows that the restrictions of \( X^{\omega _{0}} \)
to \( \gamma _{1},\cdots ,\gamma _{n} \) are independent of the choice
of volume form. It is apparent from the definition of the modular
vector field that it is unchanged if the orientation of the surface
is reversed.

Suppose that \( \pi \in \Pi (\Sigma ) \) vanishes linearly on a curve
\( \gamma  \). On a small neighborhood of \( \gamma  \), let \( \theta  \)
be the coordinate along the flow of the modular vector field \( X^{\omega _{0}} \)
with respect to \( \omega _{0} \) such that \( X^{\omega _{0}}=\partial _{\theta } \).
Since \( \pi  \) vanishes linearly on \( \gamma  \), there exists
an annular coordinate neighborhood \( (U,z,\theta ) \) of the curve
\( \gamma  \) such that\begin{eqnarray}
 &  & U=\{(z,\theta )|\, |z|<R,\, \theta \in [0,\, 2\pi ]\}, \qquad R>0, \label{cylindrical_U} \\
 &  & \gamma =\{(z,\theta )|\, z=0\},\\
 &  & \omega _{0}|_{U}=dz\wedge d\theta ,\label{omega_0_on_U} \\
 &  & \pi |_{U}=cz\partial _{z}\wedge \partial _{\theta },\quad c>0.\label{pi_on_U} 
\end{eqnarray}
Using this coordinates, it is easy to verify the following

\begin{claim}
\label{claim:orientation_modular_v_field} The restriction of a modular
vector field to a zero curve \( \gamma \in Z(\pi ) \) on which
the Poisson structure vanishes linearly is positive with respect
to the orientation on \( \gamma  \) defined by \( \pi  \) (see Lemma
 \ref{claim:orientation_curves} for the definition of this orientation.)
\end{claim}
\begin{defn}
(see also \cite{Roytenberg}) For a Poisson structure \( \pi \in \Pi (\Sigma ) \)
vanishing linearly on a curve \( \gamma \in Z(\pi ) \) define the
\emph{modular period} \emph{of} \( \pi  \) \emph{around} \( \gamma  \)
to be\[
T_{\gamma }(\pi )\doteq \textrm{period of }\, {X^{\omega _{0}}}|_{\gamma }, \]
where \( X^{\omega _{0}} \) is the modular vector field of \( \pi  \)
with respect to a volume form \( \omega _{0} \). Since \( {X^{\omega _{0}}}|_{\gamma } \)
is independent of choice of \( \omega _{0} \), the modular period
is well-defined. 
\end{defn}
Using the coordinate neighborhood \( (U,z,\theta ) \) of the curve
\( \gamma  \), we obtain\begin{equation}
\label{period_about_gamma}
T_{\gamma }(\pi )=\frac{2\pi }{c},
\end{equation}
where \( c>0 \) is as in (\ref{pi_on_U}.) 

It turns out that the modular period of the Poisson structure (\ref{pi_on_U})
on an annulus \( U \) is the only invariant under Poisson isomorphisms:

\begin{lem}
\label{modular_period_one_curve}Let \( U(R)=\{(z,\theta )|\, |z|<R,\, \theta \in [0,2\pi ]\} \)
and \( U'(R')=\{(z',\theta ')|\, |z'|<R',\, \theta '\in [0,2\pi ]\} \)
be open annuli with the orientations induced by the symplectic forms
\( \omega _{0}=dz\wedge d\theta  \) and \( \omega _{0}'=dz'\wedge d\theta ' \)
respectively. Let \( \pi =cz\partial _{z}\wedge \partial _{\theta },\, c>0 \)
and \( \pi '=c'z'\partial _{z'}\wedge \partial _{\theta '},\, c'>0 \)
be Poisson structures on \( U(R) \) and \( U'(R') \) for which the
modular periods around the zero curves \( \gamma =\{(z,\theta )|z=0\} \)
and \( \gamma '=\{(z',\theta ')|\, z'=0\} \) are equal, \( T_{\gamma }(\pi )=T_{\gamma '}(\pi ') \).
Then there is an orientation-preserving Poisson isomorphism \( \Phi :\, (U(R),\, \pi )\to (U'(R'),\, \pi ') \). 
\end{lem}
\begin{proof}
Since the modular periods are equal, we have  \( c=c' \). The map \( \Phi :\, (U(R),\pi )\to (U'(R'),\pi ') \)
given by\[
\Phi (z,\theta )=\left( \frac{R'}{R}z,\, \theta \right) \]
is a Poisson isomorphism since \( \frac{R'}{R}z\cdot \frac{R}{R'}\partial _{z}\wedge \partial _{\theta }=z\partial _{z}\wedge \partial _{\theta } \).
It is easy to see that \( \Phi  \) preserves the orientation.
\end{proof}
The fact that this Poisson isomorphism allows to change the radius
of an annuli will be used later in the proof of the classification
theorem.

\subsection{\label{section:volume}The regularized Liouville volume invariant}

To classify the topologically stable Poisson structures \( \mathcal{G}(\Sigma ) \)
up to orientation-preserving Poisson isomorphisms, we need to introduce
one more invariant.

Let \( \omega =(\tilde{\pi }^{-1}\otimes \tilde{\pi }^{-1})(\pi ) \). This form is smooth only on 
\( \Sigma \setminus Z(\pi ) \) and defines there a symplectic structure.  The
symplectic volume of each of the \( 2 \)-dimensional symplectic leaves
is infinite because the form \( \omega  \) blows up on the curves
\( \gamma _{1},\cdots ,\gamma _{n}\in Z(\pi ) \). However, there
is a way to associate a certain finite volume invariant to a Poisson
structure in \( \mathfrak {G}(\Sigma ) \), given by the principal
value of the integral\[
V(\pi )=P.V.\int _{\Sigma }\omega .\]
More precisely, let \( h\in C^{\infty }(\Sigma ) \) be a function
vanishing linearly on \( \gamma _{1},\cdots ,\gamma _{n} \) and not
zero elsewhere. Let \( \mathcal{L} \) be the set of \( 2 \)-dimensional
symplectic leaves of \( \pi  \). For \( L\in \mathcal{L} \) the
boundary \( \partial L \) is a union of curves \( \gamma _{i_{1}},\dots ,\gamma _{i_{k}}\in Z(\pi ) \).
(Note that a leaf \( L \) cannot approach the same curve from both
sides.) The function \( h \) has constant sign on each of the leaves
\( L\in \mathcal{L} \). 
 Define \[
V^{\varepsilon }_{h}(\pi )\doteq \int _{|h|>\varepsilon }\omega =\sum _{L\in \mathcal{L}}\int _{L\cap h^{-1}((-\infty ,-\varepsilon )\cup (\varepsilon ,\infty )) }\omega .\]

\begin{thm}
The limit \( V(\pi )\doteq \lim _{\varepsilon \to 0}V^{\varepsilon }_{h}(\pi ) \)
exists and is independent of the choice of function \( h \).
\end{thm}
\begin{proof}
To show that the limit in the definition of \( V(\pi ) \) exists
and is well-defined, it is enough to argue locally, in a neighborhood
of the zero set of the Poisson structure. 

For \( i=1,\dots ,n \), let \( U_{i}=\{(z_{i},\theta _{i})|\, |z_{i}|<R_{i},\, \theta _{i}\in [0,2\pi ]\} \)
be annular coordinate neighborhoods of curves \( \gamma _{i} \) such
that the restriction of \( \pi  \) on \( U_{i} \) is given by \( \pi _{|U_{i}}=c_{i}z_{i}\partial _{z_{i}}\wedge \partial _{\theta _{i}},\, c_{i}>0 \)
and \( U_{i}\cap Z(\pi )=\gamma _{i} \). Let \( \mathcal{U}=\bigsqcup _{i=1}^{n}U_{i} \).
Let \( h \) be a function vanishing linearly on the curves \( \gamma _{1},\cdots ,\gamma _{n} \)
and not zero elsewhere. Let \( z_{i}^{\pm } \) be functions of \( \theta_{i}  \)
such that \( h(\theta_{i} ,z^{\pm }_{i})=\pm \varepsilon  \). It suffices
to show that the limit\[
\lim _{\varepsilon \to 0}\textrm{P}.\textrm{V}.\int _{0}^{2\pi }\int _{z_{i}^{-}}^{z_{i}^{+}}\frac{c_{i}dz_{i}\wedge d\theta _{i}}{z_{i}}\]
is equal to zero. Indeed, this limit is equal to\[
\int _{0}^{2\pi }\left( \lim _{\varepsilon \to 0}\textrm{P}.\textrm{V}.\int _{z_{i}^{-}}^{z_{i}^{+}}\frac{c_{i}dz_{i}}{z_{i}}\right) d\theta _{i}=c_{i}\int _{0}^{2\pi }\lim _{\varepsilon \to 0}\ln \left| \frac{z_{i}^{+}}{z_{i}^{-}}\right| d\theta _{i}=0.\]

\end{proof}
Hence \( V(\pi )\in \mathbb {R} \) is a global equivalence invariant
of a Poisson structure \( \pi \in \mathfrak {G}(\Sigma ) \) on an
oriented surface which we call the \emph{regularized Liouville volume}
since in the case of a symplectic structure (i.e., \( \pi \in \mathfrak {G}_{0}(\Sigma ) \))
it is exactly the Liouville volume. If we reverse the orientation
of \( \Sigma  \), the regularized volume invariant changes sign.

\section{\label{sec:classification}A classification of topologically stable
Poisson structures}

\begin{thm}
\label{main_classification_theorem}Topologically stable Poisson structures
\( \mathfrak {G}_{n}(\Sigma ) \) on a compact  connected oriented
surface \( \Sigma  \) are completely classified (up to an orientation-preserving
Poisson isomorphism) by the following data:
\begin{enumerate}
\item The equivalence class \( [Z(\pi )]\in \mathfrak {M}_{n}(\Sigma ) \)
of the set \( Z(\pi )=\bigsqcup _{i=1}^{n}\gamma _{i} \) of zero
curves with orientations defined by \( \pi  \);
\item The modular periods around the zero curves \( \{\gamma _{i}\mapsto T_{\gamma _{i}}(\pi )|\textrm{ }i=1,\dots ,n\} \);
\item The regularized Liouville volume \( V(\pi ) \).
\end{enumerate}
In other words, two Poisson structures \( \pi ,\, \pi '\in \mathfrak {G}_{n}(\Sigma ) \)
are globally equivalent if and only if their sets of oriented zero
curves are diffeomorphism equivalent, the modular periods around the
corresponding curves are the same, and the regularized Liouville volumes
are equal.
\end{thm}
To prove this result we will need the following

\begin{lem}
\label{two_symplectic_forms}Let \( D \) be a  connected \( 2 \)-dimensional
manifold, and \( \omega _{1},\, \omega _{2} \) be two symplectic
forms on \( D \) inducing the same orientation and such that
\begin{itemize}
\item \( \omega _{1}|_{D\setminus K}=\omega _{2}|_{D\setminus K} \) for
a compact set \( K\subset D \);
\item \( \int _{D}\omega _{1}=\int _{D}\omega _{2} \);
\end{itemize}
Then there exists a symplectomorphism \( \varphi :\, (D,\omega _{1})\to (D,\omega _{2}) \)
such that \( \varphi |_{D\setminus K}=\textrm{id} \).
\begin{proof}
(Moser's trick.) Let \( \omega _{t}\doteq \omega _{1}\cdot (1-t)+\omega _{2}\cdot t \)
for \( t\in [0,1] \). Since \( \omega _{1}|_{D\setminus K}=\omega _{2}|_{D\setminus K} \),
the form \( \Delta (\omega )\doteq \omega _{2}-\omega _{1} \) is
compactly supported (\( \textrm{supp}(\Delta (\omega ))\subseteq K \).)
Since \( \int _{D}\omega _{1}=\int _{D}\omega _{2} \), the class
of \( \Delta (\omega ) \) in the second de Rham cohomology with compact
support \( H^{2}_{\textrm{c}}(D) \) is trivial. Hence \( \Delta (\omega )=d\mu  \)
for a \( 1 \)-form \( \mu \in \Omega ^{1}_{\textrm{c}}(D) \). Then
for \( v_{t}\doteq -\tilde{\omega }_{t}^{-1}(\mu ) \) we have\[
L_{v_{t}}\omega _{t}=d\imath _{v_{t}}\omega _{t}+\imath _{v_{t}}d\omega _{t}=d\imath _{v_{t}}\omega _{t}=-\Delta (\omega ). \]
Therefore, \begin{equation}
\label{v_t-equation}
L_{v_{t}}\omega _{t}+\frac{d\omega _{t}}{dt}=0.
\end{equation}
 Let \( \rho _{t} \) be the flow of the time-dependent vector field
\( v_{t} \). Since\[
\frac{d}{dt}(\rho ^{*}_{t}\omega _{t})=\rho ^{*}_{t}\left( L_{v_{t}}\omega _{t}+\frac{d\omega _{t}}{dt}\right) =0,\]
\( \rho _{t}^{*}\omega _{t}=\omega _{1} \) for all \( t\in [0,1] \).
Since \( v_{t}=0 \) outside of \( K \), it follows that \( \rho _{t}|_{D\setminus K}=\textrm{id} \).
Define \( \varphi =\rho _{1} \). Then \( \varphi _{|D\setminus K}=\textrm{id} \)
and \( \varphi ^{*}\omega _{2}=\rho ^{*}_{1}\omega _{2}=\omega _{1} \)
as desired. 
\end{proof}
\end{lem}
We now have all of the ingredients for the proof of  
Theorem~\ref{main_classification_theorem}.

\begin{proof}
(of Theorem \ref{main_classification_theorem}.) Let \( \pi  \) and
\( \pi ' \) be two topologically stable Poisson structures. Clearly,
the coincidence of all of the invariants listed in the theorem is
necessary for \( \pi  \) and \( \pi '  \) to be isomorphic. 

Conversely, assume that all the invariants are the same. Then the
zero sets \( Z(\pi )\subset \Sigma  \) and \( Z(\pi ')\subset \Sigma  \)
are diffeomorphic and, therefore, we may assume that \( Z(\pi )=Z(\pi ') \)
by replacing \( \pi ' \) with an isomorphic structure. Since both
\( \pi =f\cdot \pi _{0} \) and \( \pi '=f'\cdot \pi _{0} \) vanish
linearly on each connected component \( \gamma \subset Z(\pi )=Z(\pi ') \),
we can, by once again replacing \( \pi ' \) by a Poisson-isomorphic
structure, assume that in an  annular coordinate neighborhood \( U_i \simeq S^{1}\times (-1,1)=\{(z_i,\theta _i )|\, \theta _i\in [0, 2\pi] ,\, |z_i|<R_i\} \)
of \( \gamma _i  \) we have\begin{eqnarray*}
\pi |_{U_{i}}  & = & c_iz_i\partial _{z_{i}}\wedge \partial _{\theta },\quad c_i\neq 0,\\
\pi '|_{U_{i}} & = & c'_iz_i\partial _{z_i}\wedge \partial _{\theta },\quad c'_i\neq 0.
\end{eqnarray*}
Since the modular periods of \( \pi  \) and \( \pi ' \) around the
corresponding curves are assumed to be the same, we have \( 2\pi /|c_i|=2\pi /|c'_i| \),
which implies \( |c_i|=|c'_i| \). The fact that the orientations of the
connected components of \( Z(\pi ) \) and \( Z(\pi ') \) induced
by modular vector fields are the same implies that \( c_i=c'_i \). Hence
\( \pi  \) and \( \pi ' \) are equal in a neighborhood of \( Z(\pi ) \).
By replacing \( \pi ' \) with an isomorphic structure we may assume
that \( \pi =\pi ' \) on a neighborhood \( U  \) of \( Z(\pi ) \).

Consider the non-compact manifold \( N\doteq \Sigma \setminus Z(\pi ) \).
The desired isomorphism of \( \pi  \) and \( \pi ' \) on \( \Sigma  \)
would follow if we could argue that \( \pi  \) and \( \pi ' \) are
isomorphic on \( N \) via a diffeomorphism \( \varphi  \) which
is an identity on \( N\cap U \). Indeed, \( \varphi  \) could then
be extended to a diffeomorphism \( \bar{\varphi } \) of \( \Sigma  \)
by setting \( \varphi (p)=p \) for \( p\in Z(\pi ) \) and we would
clearly have \( \bar{\varphi }_{*}(\pi )=\pi ' \). Since \( \pi  \)
and \( \pi ' \) are nonzero on \( N \), they come from symplectic
structures \( \frac{\omega _{0}}{f} \) and \( \frac{\omega _{0}}{f'} \)
on \( N \). These \( 2 \)-forms coincide on \( N\cap U \). Moreover,
\( \int _{N}(\frac{\omega _{0}}{f}-\frac{\omega _{0}}{f'})=V(\pi )-V(\pi ')=0 \).
Hence \( \frac{\omega _{0}}{f} \) and \( \frac{\omega _{0}}{f'} \)
define the same class in the compactly-supported de Rham cohomology
of \( N \). It is now a matter of imitating the proof of Moser's
theorem to construct a time-dependent vector field supported on \( N\setminus U \)
and whose flow at time \( 1 \) carries \( \frac{\omega _{0}}{f} \)
into \( \frac{\omega _{0}}{f'} \). Taking \( \varphi  \) to be that
flow at time \( 1 \) we find that we have obtained the desired isomorphism.
\end{proof}
\begin{rem}
Given \( \pi ,\pi '\in \mathcal{G}_{n}(\Sigma ) \), one can ask if
\( (\Sigma ,\pi ) \) and \( (\Sigma ,\pi ') \) are equivalent by
an arbitrary (possibly orientation-reversing) Poisson isomorphism.
Fix an orientation-reversing diffeomorphism \( \nu :\Sigma \to \Sigma  \).
Then \( \pi  \) and \( \pi ' \) are equivalent by an orientation-reversing
diffeomorphism if and only if \( \nu _{*}\pi  \) and \( \pi ' \)
are equivalent by an orientation-preserving diffeomorphism. It is
not hard to see that \( T_{\gamma _{i}}(\pi )=T_{\nu (\gamma _{i})}(\nu _{*}\pi ) \)
for all \( \gamma _{i}\in Z(\pi ) \) and \( V(\pi )=-V(\nu _{*}\pi ) \).
Thus the question of equivalence by orientation reversing maps can
be reduced to the orientation-preserving context of Theorem \ref{main_classification_theorem}. 
\end{rem}
\begin{example}
\label{Example: not_globally_equivalent}Let \( \omega  \) and \( \omega '=-\omega  \)
be two symplectic structures on a compact oriented surface. Then \( \omega  \)
and \( \omega ' \) are Poisson isomorphic by an orientation-reversing
diffeomorphism, but not by an orientation-preserving diffeomorphism.

There are, of course, similar examples of structures with non-trivial
sets of linear degeneracy. Consider the unit \( 2 \)-sphere \( S^{2} \)
with the cylindrical polar coordinates \( (z,\theta ) \) away from
its poles. Let \( \omega _{0}=dz\wedge d\theta  \) be a symplectic
form on \( S^{2} \) with the corresponding Poisson bivector \( \pi _{0} \).
Let \( \pi ,\, \pi '\in \mathscr {G}_{2}(S^{2}) \) be the Poisson
structures given by\[
\pi =(z-a)(z-b)\partial _{z}\wedge \partial _{\theta },\quad -1<b<a<1\]
 and \( \pi '=-\pi  \). Choose \( a \) and \( b \) in such a way
that \( V(\pi )=V(\pi ')=0 \). Let \( \gamma _{1}=\{(z,\theta )|\, z=a\} \)
and \( \gamma _{2}=\{(z,\theta )|\, z=b\} \) be the zero curves of
\( \pi ,\, \pi ' \). On both \( \gamma _{1} \) and \( \gamma _{2} \)
the orientations defined by \( \pi  \) and \( \pi ' \) are opposite
to each other. Let \( L_{\textrm{top}}=\{(z,\theta )|\, a<z<1\} \),
\( L_{\textrm{middle}}=\{(z,\theta )|\, b<z<a\} \) and \( L_{\textrm{bottom}}=\{(z,\theta )|\, -1<z<b\} \)
be the \( 2 \)-dimensional leaves (common for both structures). The
structures \( \pi  \) and \( \pi ' \) can not be Poisson isomorphic
in an orientation-preserving way since such a diffeomorphism would
have to exchange the two-dimensional disks \( L_{\textrm{top }} \)and
\( L_{\textrm{bottom }} \) with the annulus \( L_{\textrm{middle}} \).
On the other hand, \( (S^{2},\pi ) \) and \( (S^{2},\pi ') \) are
clearly Poisson isomorphic by an orientation-reversing diffeomorphism
\( (z,\theta )\mapsto (z,-\theta ) \).
\end{example}

\section{\label{sec:cohomology}Poisson cohomology of topologically stable
Poisson structures }

In this section we compute the Poisson cohomology of a given topologically
stable Poisson structure on a compact connected oriented surface and
describe its relation to the infinitesimal deformations and the classifying
invariants introduced above. (For generalities on Poisson cohomology
see, e.g., \cite{Vaisman}.) 

First, recall the following

\begin{lem}
\label{H_of_annulus}(e.g., Roytenberg \cite{Roytenberg}) The Poisson
cohomology of an annular neighborhood \( U=\{(z,\theta )|\, |z|<R,\, \theta \in [0,2\pi ]\} \)
of a zero  curve \( \gamma  \) on which \( \pi _{|U}=z\partial _{z}\wedge \partial _{\theta } \)
is given by \begin{eqnarray*}
 &  & H^{0}_{\pi }(U,\pi _{|U})=\textrm{span}\langle 1\rangle =\mathbb {R},\\
 &  & H^{1}_{\pi }(U,\pi _{|U})=\textrm{span}\langle \partial _{\theta },\, z\partial _{z}\rangle =\mathbb {R}^{2},\\
 &  & H^{2}_{\pi }(U,\pi _{|U})=\textrm{span}\langle z\partial _{z}\wedge \partial _{\theta }\rangle =\mathbb {R}.
\end{eqnarray*}

\end{lem}
Thus, \( H^{0}_{\pi }(U) \) is generated by constant functions. The
first cohomology \( H^{1}_{\pi }(U) \) is generated by the modular
class (\( \partial _{\theta } \) is the modular vector field of \( \pi _{|U} \)
with respect to \( \omega _{0}=dz\wedge d\theta  \)) and the image
of the first de Rham cohomology class of \( U \) (spanned by \( d\theta  \))
under the homomorphism \( \tilde{\pi }:\, H^{1}(U)\to H^{1}_{\pi }(U) \),
which is injective in this case. The second cohomology is generated
by \( \pi _{|U} \).

Let \( \pi \in \mathfrak {G}_{n}(\Sigma ) \) be a topologically stable
Poisson structure on \( \Sigma  \). Since a Casimir function on \( \Sigma  \)
must be constant on all connected components of \( \Sigma \setminus Z(\pi ) \),
by continuity it must be constant everywhere. Hence \( H^{0}_{\pi }(\Sigma )=\mathbb {R}=\textrm{span}\langle 1\rangle  \). 

We will (inductively) use the Mayer-Vietoris sequence of Poisson cohomology
(see, e.g., \cite{Vaisman}) to compute \( H^{1}_{\pi }(\Sigma ) \)
and \( H^{2}_{\pi }(\Sigma ) \). 

Let \( U_{i} \) be an annular neighborhood of the curve \( \gamma _{i}\in Z(\pi ) \)
such that \( U_{i}\cap Z(\pi )=\gamma _{i} \). Let \( V_{0}\doteq \Sigma  \)
and define inductively \( V_{i}\doteq V_{i-1}\setminus \gamma _{i} \)
for \( i=1,\ldots ,n \). Consider the cover of \( V_{i-1} \) by
open sets \( U_{i} \) and \( V_{i} \). Consider the exact Mayer-Vietoris
sequence of Poisson cohomology associated to this cover:\[
\begin{array}{cccccccc}
0 & \to  & H^{0}_{\pi }(V_{i-1}) & \stackrel{\alpha ^{0}_{i}}{\to } & H^{0}_{\pi }(U_{i})\oplus H^{0}_{\pi }(V_{i}) & \stackrel{\beta ^{0}_{i}}{\to } & H^{0}_{\pi }(U_{i}\cap V_{i}) & \stackrel{\delta ^{0}_{i}}{\to }\\
 & \to  & H^{1}_{\pi }(V_{i-1}) & \stackrel{\alpha ^{1}_{i}}{\to } & H^{1}_{\pi }(U_{i})\oplus H^{1}_{\pi }(V_{i}) & \stackrel{\beta ^{1}_{i}}{\to } & H^{1}_{\pi }(U_{i}\cap V_{i}) & \stackrel{\delta ^{1}_{i}}{\to }\\
 & \to  & H^{2}_{\pi }(V_{i-1}) & \stackrel{\alpha ^{2}_{i}}{\to } & H^{2}_{\pi }(U_{i})\oplus H^{2}_{\pi }(V_{i}) & \stackrel{\beta ^{2}_{i}}{\to } & H^{2}_{\pi }(U_{i}\cap V_{i}) & \to 0.
\end{array}\]
By exactness, \( H^{1}_{\pi }(V_{i-1})\simeq \textrm{Im}\left( \delta ^{0}_{i}\right) \oplus \textrm{ker}\beta ^{1}_{i} \),
where\[
\beta _{i}^{1}\left( [\chi ]_{U_{i}},\, [\nu ]_{V_{i}}\right) =[\chi -\nu ]_{U_{i}\cap V_{i}}, \quad \chi \in \mathfrak {X}^{1}_{\pi }(U_{i}),\, \nu \in \mathfrak {X}^{1}_{\pi }(V_{i}),\, d_{\pi }\chi = d_{\pi }\nu =0,\]
and \( [X]_{W} \) denotes the class of the (Poisson) vector field
\( X|_{W} \) in \( H^{1}_{\pi }(W) \). 

By Lemma \ref{H_of_annulus}, \( H^{1}_{\pi }(U_{i})\simeq \tilde{\pi }(H^{1}(U_{i}))\oplus \textrm{span}\langle \partial _{\theta _{i}}\rangle \simeq \mathbb {R}^{2} \).
Since \( U_{i}\cap V_{i} \) is a union of two symplectic annuli,
\( H^{1}_{\pi }(U_{i}\cap V_{i})=\tilde{\pi }(H^{1}(U_{i}\cap V_{i}))=\mathbb {R}^{2} \).
Therefore, \begin{equation}
\label{H_Poisson_inductive}
H^{1}_{\pi }(V_{i-1})\simeq \delta ^{0}_{i}(H^{0}_{\pi }(U_{i}\cap V_{i}))\oplus \textrm{span}\langle \partial _{\theta _{i}}\rangle \oplus \textrm{ker}\left( \beta ^{1}_{i}|_{\tilde{\pi }(H^{1}_{\textrm{ }}(U_{i}))\oplus H^{1}_{\pi }(V_{i})}\right) .
\end{equation}
Consider also the long exact sequence in de Rham cohomology associated
to the same cover\[
\begin{array}{cccccccc}
0 & \to  & H^{0}(V_{i-1}) & \stackrel{a^{0}_{i}}{\to } & H^{0}(U_{i})\oplus H^{0}(V_{i}) & \stackrel{b^{0}_{i}}{\to } & H^{0}(U_{i}\cap V_{i}) & \stackrel{d^{0}_{i}}{\to }\\
 & \to  & H^{1}(V_{i-1}) & \stackrel{a^{1}_{i}}{\to } & H^{1}(U_{i})\oplus H^{1}(V_{i}) & \stackrel{b^{1}_{i}}{\to } & H^{1}(U_{i}\cap V_{i}) & \stackrel{d^{1}_{i}}{\to }\\
 & \to  & H^{2}(V_{i-1}) & \stackrel{a^{2}_{i}}{\to } & H^{2}(U_{i})\oplus H^{2}(V_{i}) & \stackrel{b^{2}_{i}}{\to } & H^{2}(U_{i}\cap V_{i}) & \to 0.
\end{array}\]
By exactness, we have \( H^{1}(V_{i-1})\simeq \textrm{Im}\left( d^{0}_{i}\right) \oplus \textrm{ker}b^{1}_{i} \).
Since \( V_{n}=\Sigma \setminus Z(\pi ) \) is symplectic, \( H^{*}_{\pi }(V_{n})=\tilde{\pi }(H^{*}(V_{n})) \).
This together with \( H^{0}_{\pi }(U_{n})\simeq H^{0}(U_{n}) \),
\( H^{0}_{\pi }(U_{n}\cap V_{n})=\tilde{\pi }(H^{0}(U_{n}\cap V_{n})) \)
implies \( \textrm{Im}(\delta ^{0}_{n})=\tilde{\pi }(\textrm{Im}(d^{0}_{n})) \)
and, therefore,\[
\textrm{ker}\left( \beta ^{1}_{n}|_{\tilde{\pi }(H^{1}_{\textrm{ }}(U_{n})\oplus H^{1}_{\pi }(V_{n}))}\right) =\tilde{\pi }(\textrm{ker}(b^{1}_{n})).\]
Hence, from (\ref{H_Poisson_inductive}) it follows that \[
H^{1}_{\pi }(V_{n-1})\simeq \tilde{\pi }(\textrm{Im}(d_{n-1}))\oplus \textrm{span}\langle \partial _{\theta _{n}}\rangle \oplus \tilde{\pi }(H^{1}(V_{n-1})).\]
For \( i=n-2 \), we have\begin{multline*}
H^{1}_{\pi }(V_{n-2})\simeq \textrm{span}\langle \partial _{\theta _{n-1}}\rangle \oplus \textrm{ker}\left( \beta ^{1}_{n-1}|_{\tilde{\pi }(H^{1}_{\textrm{ }}(U_{n-1})\oplus \left( \tilde{\pi }(H^{1}(V_{n-1}))\oplus \tilde{\pi }(\textrm{Im}(d^{0}_{i}))\right) }\right) = \\
=\textrm{span}\langle \partial _{\theta _{n-1}}\rangle \oplus \textrm{span}\langle \partial _{\theta _{n}}\rangle \oplus \tilde{\pi }(H^{1}(V_{n-2})). 
\end{multline*}
Working inductively (from \( i=n-1 \) to \( i=0 \)), we obtain\[
H^{1}_{\pi }(\Sigma )\simeq \mathbb {R}^{2g+n}=H^{1}(\Sigma )\oplus \bigoplus _{i=1}^{n}\textrm{span}\langle \partial _{\theta _{i}}\rangle ,\]
where \( g \) is the genus of the surface \( \Sigma  \). 

For the second Poisson cohomology, we have \( H^{2}_{\pi }(U_{i}\cap V_{i})=0 \),
\( H^{2}_{\pi }(U_{i})=\mathbb {R} \) for all \( i=1,\ldots ,n \)
and \( H^{2}_{\pi }(V_{n})=0 \). Therefore, applying the Mayer-Vietoris
sequence inductively, we obtain\[
H^{2}_{\pi }(V_{i-1})\simeq \mathbb {R}\oplus H^{2}_{\pi }(U_{i})\oplus H^{2}_{\pi }(V_{i}),\]
 where \( \mathbb {R} \) is generated by the image of \( (r_{i}\partial _{r_{i}},\, -r_{i}\partial _{r_{i}})\in H^{1}_{\pi }(U_{i}\cap V_{i}) \)
under the connecting homomorphism \( \delta ^{0}_{i} \). One can
show (similarly to Lemma 4.3.1 in \cite{Roytenberg}) that the class
of the standard non-degenerate Poisson structure \( \pi _{0}|_{V_{i-1}} \)
is not trivial in \( H^{2}_{\pi }(V_{i-1}) \). Therefore,\[
H^{2}_{\pi }(\Sigma )\simeq \mathbb {R}^{n+1}=\bigoplus _{i=1}^{n}H^{2}_{\pi }(U_{i})\oplus \mathbb {R},\]
where the first \( n \) generators are the Poisson structures of
the form\begin{equation}
\label{pi_i}
\pi _{i}\doteq \pi \cdot \{\textrm{bump function around the curve }\gamma _{i}\},\quad i=1,\ldots ,n,
\end{equation}
and the last generator is the standard non-degenerate Poisson structure
\( \pi _{0} \) on \( \Sigma  \). Therefore, we have proved the following

\begin{thm}
Let \( \pi \in \mathfrak {G}_{n}(\Sigma ) \) be a topologically stable
Poisson structure on a compact connected oriented surface \( \Sigma  \)
of genus \( g \). The Poisson cohomology of \( \pi  \) is given
by\begin{eqnarray*}
 &  & H^{0}_{\pi }(\Sigma ,\pi )=\textrm{span}\langle 1\rangle =\mathbb {R},\\
 &  & H^{1}_{\pi }(\Sigma ,\pi )=\textrm{span}\langle X^{\omega _{0}}(\pi _{1}),\cdots ,X^{\omega _{0}}(\pi _{n})\rangle \oplus \tilde{\pi }(H^{1}_{\textrm{deRham}}(\Sigma ))=\mathbb {R}^{n+2g},\\
 &  & H^{2}_{\pi }(\Sigma ,\pi )=\textrm{span}\langle \pi _{0};\, \pi _{1},\cdots ,\pi _{n}\rangle =\mathbb {R}^{n+1},
\end{eqnarray*}
where \( \pi _{0} \) is a non-degenerate Poisson structure on \( \Sigma  \)
, \( \pi _{i}\in \mathfrak {G}_{1}(\Sigma ),\, i=1,\ldots ,n \) is
a Poisson structure vanishing linearly on \( \gamma _{i}\in Z(\pi ) \)
and identically zero outside of a neighborhood of \( \gamma _{i} \); and 
\( X^{\omega _{0}}(\pi _{i}) \) is the modular vector field of \( \pi _{i} \)
with respect to the standard symplectic form \( \omega _{0} \) on
\( \Sigma  \).
\end{thm}
Notice that the dimensions of the cohomology spaces depend only on
the number of the zero curves and not on their positions. In particular,
the Poisson cohomology as a vector space does not depend on the homology
classes of the zero curves of the structure. Recall (see, e.g., \cite{Vaisman})
that the Poisson cohomology space \( H^{*}_{\pi }(P) \) has the structure
of an associative graded commutative algebra induced by the operation
of wedge multiplication of multivector fields. A direct computation
verifies the following

\begin{prop}
The wedge product on the cohomology space \( H^{*}_{\pi }(\Sigma ,\pi ) \)
of a topologically stable  Poisson structure on \( \Sigma  \) is
determined by\begin{eqnarray*}
 &  & [1]\wedge [\chi ]=[\chi ],\quad [\chi ]\in H^{*}_{\pi }(\Sigma ,\pi ),\\
 &  & [X^{\omega _{0}}(\pi _{i})]\wedge [X^{\omega _{0}}(\pi _{j})]=0,\quad i,j=1,\dots ,n,\\
 &  & [X^{\omega _{0}}(\pi _{i})]\wedge [\tilde{\pi }(\alpha )]=[\alpha (X^{\omega _{0}}(\pi _{i}))\cdot \pi _{i}]=\left( \frac{1}{T_{\gamma _{i}}(\pi )}\int _{\gamma _{i}}\alpha \right) \cdot [\pi _{i}],\\
 &  & [\tilde{\pi }(\alpha )]\wedge [\tilde{\pi }(\alpha ')]=\tilde{\pi }(\bar{\alpha }\wedge \overline{\alpha '}),\\
 &  & [\chi ]\wedge H^{2}_{\pi }(\Sigma ,\pi )=0,\quad [\chi ]\in H^{1}_{\pi }(\Sigma )\oplus H^{2}_{\pi }(\Sigma ),
\end{eqnarray*}
where \( \alpha ,\alpha '\in \Omega ^{1}(\Sigma ),\textrm{ }d\alpha =d\alpha '=0 \). 
\end{prop}
(Here bar denotes the class of its argument in the de Rham cohomology
and the brackets \( [\, \, ] \) denote the class in the Poisson cohomology.)
We should mention that the wedge product \( \wedge  \) in de Rham
cohomology is dual to the intersection product in homology. 

This computation allows one to compute the number of zero curves \( \gamma _{k} \),
which determine non-zero homology classes. To see this, we note that
\( [\pi _{k}]\in H^{1}_{\pi }(\Sigma )\wedge H^{1}_{\pi }(\Sigma ) \)
iff there exists a \( 1 \)-form \( \alpha  \) such that \( \int _{\gamma _{k}}\alpha \neq 0 \),
i.e., \( \gamma _{k} \) is non-zero in homology. If \( \Sigma  \)
is not a sphere, \( H^{1}_{\textrm{deRham}}(\Sigma ) \) is non-zero.
Since the intersection form on \( H^{1}_{\textrm{deRham}}(\Sigma ) \)
is non-degenerate (implementing Poincare duality), it follows that
\( H^{1}_{\textrm{deRham}}(\Sigma )\wedge H^{1}_{\textrm{deRham}}(\Sigma )\neq 0 \).
Thus in the case that \( \Sigma  \) is not a sphere, \( H^{1}_{\pi }(\Sigma )\wedge H^{1}_{\pi }(\Sigma ) \)
has the set\[
\{\pi _{0}\}\cup \{\pi _{k}:\gamma _{k}\textrm{ is nontrivial in homology}\}\]
as a basis and so the number of curves \( \gamma _{k} \), which are
non-trivial in homology, is just \( \dim (H^{1}_{\pi }(\Sigma )\wedge H^{1}_{\pi }(\Sigma ))-1 \).
In the case that \( \Sigma  \) is a sphere, all \( \gamma _{k} \)
are of course topologically trivial. 

We can now interpret the generators of \( H^{2}_{\pi }(\Sigma ,\pi ) \)
as infinitesimal deformations of the Poisson structure \( \pi  \)
and find out how these deformations affect the classifying invariants.

\begin{cor}
Let \( \pi \in \mathfrak {G}_{n}(\Sigma ) \), \( Z(\pi )=\bigsqcup _{i=1}^{n}\gamma _{i} \).
The following \( n+1 \) one-parameter families of infinitesimal deformations
form a basis of \( H^{2}_{\pi }(\Sigma ,\pi ) \)

(1) \( \pi \mapsto \pi +\varepsilon \cdot \pi _{0} \);

(2) \( \pi \mapsto \pi +\varepsilon \cdot \pi _{i} \), \( i=1,\ldots ,n \).

Each of these deformations changes exactly one of the classifying
invariants of the Poisson structure: \( \pi \mapsto \pi +\varepsilon \cdot \pi _{0} \)
changes the regularized Liouville volume and \( \pi \mapsto \pi +\varepsilon \cdot \pi _{i} \)
changes the modular period around the curve \( \gamma _{i} \) for
each \( i=1,\dots ,n \). 
\end{cor}

\section{\label{sec:sphere}Example: topologically stable Poisson structures
on the sphere}

It would be interesting to describe the moduli space of the space
of topologically stable Poisson structures on an oriented surface
up to orientation-preserving diffeomorphisms. The first step would
be the description of the moduli space \( \mathfrak {M}_{n} \) of
\( n \) disjoint oriented curves on \( \Sigma  \). However, this
problem is already quite difficult for a general surface. Here we
will  consider the simplest example of topologically stable
Poisson structures on the sphere. 

Let \( (\gamma _{1},\cdots ,\gamma _{n}) \) be a set of disjoint
curves on \( S^{2} \). Let \( \mathcal{S}_{1},\cdots ,\mathcal{S}_{n+1} \)
be the connected components of \( S^{2}\setminus (\gamma _{1},\cdots ,\gamma _{n}) \).
To the configuration of curves \( (\gamma _{1},\cdots ,\gamma _{n}) \)
we associate a graph \( \Gamma (\gamma _{1},\cdots ,\gamma _{n}) \)
in the following way. The vertices \( v_{1},\dots ,v_{n+1} \) of
the graph correspond to the connected components \( \mathcal{S}_{1},\dots ,\mathcal{S}_{n+1} \).
Two vertices \( v_{i} \) and \( v_{j} \) are connected by an edge
\( e_{k} \) iff \( \gamma _{k} \) is the common bounding curve of
the regions \( \mathcal{S}_{i} \) and \( \mathcal{S}_{j} \).

\begin{claim}
For a set of disjoint curves \( \gamma _{1},\dots ,\gamma _{n} \)
on \( S^{2} \) the graph \( \Gamma (\gamma _{1},\dots ,\gamma _{n}) \)
is a tree.
\end{claim}
\begin{proof}
Let \( e_{i}\in E(\Gamma (\gamma _{1},\dots ,\gamma _{n})) \) be
an edge of the graph corresponding to the curve \( \gamma _{i} \).
Since \( S^{2}\setminus \gamma _{i} \) is a union of two open sets,
it follows that \( \Gamma (\gamma _{1},\dots ,\gamma _{n})\setminus e_{i} \)
(i.e., the graph \( \Gamma (\gamma _{1},\dots ,\gamma _{n}) \) with
the edge \( e_{i} \) removed) is a union of two disjoint graphs.
Since this is true for any \( e_{i},\, i=1,\dots ,n \), the graph
is a tree. 
\end{proof}
Choose an orientation on \( S^{2} \) and a symplectic form \( \omega _{0} \)
(with the Poisson bivector \( \pi _{0} \)) which induces this orientation.
Let \( \pi \in \mathcal{G}_{n}(\Sigma ) \) be a topologically stable
Poisson structure. The function \( f=\pi /\pi _{0} \) has constant
signs on the \( 2 \)-dimensional symplectic leaves. Let \( \Gamma (\gamma _{1},\dots ,\gamma _{n}) \)
be the tree associated to the zero curves \( \gamma _{1},\dots ,\gamma _{n} \)
of \( \pi  \) as described above. Assign to each vertex \( v_{i} \)
a sign (plus or minus) equal to the sign \emph{}of the function \( f \)
on the corresponding symplectic leaf \( \mathcal{S}_{i} \) of \( \pi  \).
The properties of \( \pi  \) imply that for any edge of this tree
its ends are assigned the opposite signs. We will call the tree associated
to the zero curves of \( \pi  \) with signs associated to its vertices
the \emph{signed tree \( \Gamma (\pi ) \)} of the Poisson structure
\emph{\( \pi  \).}

Consider the map \( T_{e}(\pi ):\, E(\Gamma (\pi ))\to \mathbb {R}^{+} \)
which for each edge \( e\in E(\Gamma (\pi )) \) gives a period \( T_{e}(\pi ) \)
of a modular vector field of \( \pi  \) around the zero curve corresponding
to this edge. The classification Theorem \ref{main_classification_theorem}
implies

\begin{thm}
The topologically stable Poisson structures \( \pi \in \mathfrak {G}_{n}(S^{2}) \)
on the sphere are completely classified (up to an orientation-preserving
Poisson isomorphism) by the signed tree \( \Gamma (\pi ) \), the
map \( e\mapsto T_{e}(\pi ) \), \( e\in E(\Gamma (\pi )) \) and
the regularized Liouville volume \( V(\pi ) \). In other words, \( \pi _{1},\, \pi _{2}\in \mathfrak {G}_{n}(S^{2}) \)
are globally equivalent if and only if the corresponding \( \Gamma (\pi _{i}) \),
\( \{e\mapsto T_{e}(\pi _{i})\}, \) \( V(\pi _{i}) \) are the same
(up to automorphisms of signed trees with positive numbers attached
to their edges.) 
\end{thm}
The moduli space of generic Poisson structures in \( \mathfrak {G}_{n}(S^{2}) \)
up to Poisson isomorphisms is\[
\mathfrak {G}_{n}(S^{2})/(\textrm{Poisson isomorphisms})\simeq \left( \bigsqcup _{\Gamma _{n+1}}(\mathbb {R}^{+})^{n}/\textrm{Aut}(\Gamma ,T_{e})\right) \times \mathbb {R},\]
 where \( \textrm{Aut}(\Gamma ,T_{e}) \) is the automorphism group
of the signed tree \( \Gamma  \) with \( n+1 \) vertices and with
positive numbers \( T_{e} \) attached to its edges \( e\in E(\Gamma (\pi )) \).
The moduli space has dimension \( n+1 \) and is coordinatized by
\( \{T_{e}(\pi )|\, e\in E(\Gamma (\pi ))\} \) and \( V(\pi ) \).

A particular case of topologically stable Poisson structures on \( S^{2} \),
the \( SU(2) \)-covariant structures vanishing on a circle on \( S^{2} \),
were considered by D.Roytenberg in \cite{Roytenberg}. In the coordinates
\( (z,\theta ) \) on the unit sphere these structures are given by\[
\pi _{b}=a(z-b)\partial _{z}\wedge \partial _{\theta }\qquad \textrm{for }|b|<1,\, a>0.\]
The modular period around the zero curve (a {}``horizontal'' circle
\( \gamma =\{(z,\theta )|\, z=b\} \)) and the regularized Liouville
volume are given by\begin{eqnarray*}
 &  & T_{\gamma =\{(z,\theta )|\, z=b\}}(\pi )=\frac{2\pi }{a},\\
 &  & V(\pi )=\frac{2\pi }{a}\ln \frac{1+b}{1-b}.
\end{eqnarray*}
Note that for a non-degenerate Poisson structure \( \pi _{b}=a(z-b)\partial _{z}\wedge \partial _{\theta },\, |b|>1 \)
the total Liouville volume is given by the same formula, \( V(\pi )=\frac{2\pi }{a}\ln \left| \frac{1+b}{1-b}\right|  \). 

\begin{cor}
Let \( T\in \mathbb {R}^{+} \) and \( V\in \mathbb {R} \). A Poisson
structure \( \pi \in \mathfrak {G}_{1}(S^{2}) \) with the modular
period \( T \) and the regularized total volume \( V \) is globally
equivalent to the Poisson structure which in coordinates \( (z,\theta ) \)
on \( S^{2} \) is given by\[
\pi (T,V)=\frac{2\pi }{T}\left( z-\frac{e^{V/T}-1}{e^{V/T}+1}\right) \partial _{z}\wedge \partial _{\theta }\]
 and vanishes linearly on the circle \( z=\frac{e^{V/T}-1}{e^{V/T}+1} \). 
\end{cor}
Utilizing Poisson cohomology, D. Roytenberg \cite[Corollary 4.3.3, 4.3.4]{Roytenberg}
has previously obtained that the structures \( \pi _{b} \), \( -1<b<1 \)
are non-trivial infinitesimal deformations of each other. Similarly,
he proved that for each \( b \), \( \pi _{b} \) admits no infinitesimal
rescalings. Using Theorem \ref{main_classification_theorem}, we get
the following improvement of his results:

\begin{cor}
(a) The Poisson structures \( \pi _{b} \) and \( \pi _{b'} \) are
globally equivalent iff \( b=b' \). \\
(b) For \( \alpha \in \mathbb {R}\setminus \{0\} \), the Poisson
structures \( \pi _{b} \) and \( \alpha \pi _{b} \) are equivalent
via an orientation-preserving Poisson isomorphism (respectively, arbitrary
Poisson isomorphism) if and only if \( \alpha =1 \) (respectively,
\( |\alpha |=1 \).) In particular, \( \pi _{b} \) admits no rescalings.
\end{cor}
\bibliographystyle{amsplain}

\def\cprime{$'$} \def\cprime{$'$}
\providecommand{\bysame}{\leavevmode\hbox to3em{\hrulefill}\thinspace}

\end{document}